\documentclass[11pt]{article}
\usepackage{amsfonts}
\usepackage{amssymb,amsfonts,amsmath,amsthm}
\usepackage{epsfig}
\newtheorem{thm}{Theorem}

\newtheorem{lem}[thm]{Lemma}
\newtheorem{prop}[thm]{Proposition}

\makeatletter \@addtoreset{equation}{section}

\def\pf{\noindent {\it Proof.\ }}
\def\qed{\hfill \rule{4pt}{7pt}}
\def \S {{\mathfrak{S}}}
\def \sg {{\sigma}}
\def \E {{\mathcal{E}}}

\def \G {{\mathcal {G}}}
\begin{document}
\begin{center}
{\large\bf A Combinatorial Proof of the Enumeration of Alternating Permutations with Given Peak Set}\\
[7pt]
Alina F.Y. Zhao\\
[5pt]

School of Mathematical Science\\Nanjing Normal University, Nanjing
210046, P.R. China \\

alinazhao@njnu.edu.cn

\end{center}
\vskip 3mm

\begin{abstract}
Using the correspondence between a cycle up-down permutation and a pair of matchings, we give a combinatorial proof of the enumeration of alternating permutations according to the given peak set.
\end{abstract}



\section{Introduction}
Let $\S_n$ denote the symmetric group of all permutations of $[n]:= \{1, 2, \ldots, n\}$. An alternating permutation on $[n]$ is defined to be a permutation
$\sg=\sg_1 \sg_2 \cdots \sg_n$ $\in \S_n$ satisfying $ \sg_1>\sg_2<\sg_3>\sg_4<\cdots$, etc., in an
alternating way. Similarly, $\sg$ is reverse alternating if $ \sg_1<\sg_2>\sg_3<\sg_4>\cdots$, which is also referred as an up-down permutation.
Denote by $\E_n$ the set of alternating permutations on $[n]$, and further let $E_n=|\E_n|$. Note that $E_n$ is called Euler number, and was shown by Andr\'{e} \cite{Andre,Andre1881} to satisfy\[
\sum_{n\geq 0}E_n \frac{x^n}{n!}=\sec x+\tan x.
\]

The reverse map $\sg \mapsto \sg^{r}$ defined by $\sg^{r}_i=\sg_{n+1-i}$ on $\S_n$ shows that $E_n$ is also the number of up-down permutations in $\S_n$. Recently, Elizalde and Deutsch \cite{Eli} introduced the concept of cycle up-down permutations. A cycle is said to be up-down if, when written in standard cycle form, say $(a_1,a_2, a_3,\ldots)$, one has $a_1<a_2> a_3<a_4>\cdots$, and a permutation $\sg$ is a cycle up-down permutation if it is a product of up-down cycles. They prove both bijectively and analytically that
\begin{prop}[\cite{Eli}, Lemma 2.2]
The number of cycle up-down permutations of $[2k]$ all of whose cycles are even is
$E_{2k}$.
\end{prop}
For out purpose, let us briefly recall the bijection $\tau$ developed in \cite{Eli} to prove the above proposition. Given $\sg=\sg_1 \sg_2 \cdots \sg_{2k}$ $\in \E_{2k}$, let $\sg_{i_1} > \sg_{i_2} >\cdots >\sg_{i_m}$ be its left to right minima, the corresponding cycle up-down permutation $\tau(\sg)$ with only even cycles is defined by
\[ \tau(\sg) = (\sg_{i_1},\ldots,\sg_{i_2-1})(\sg_{i_2},\ldots,\sg_{i_3-1}) \cdots (\sg_{i_m}, \ldots,\sg_{2k}).
\]

The element $\sg_i(1\leq i \leq n)$ is called a peak if $\sg_{i-1}<\sg_i>\sg_{i+1}$, where we set $\sg_0=0$ and $\sg_{n+1}=0$, and the peak set of $\sg$ are the elements of peaks in $\sg$. For other definitions of peaks see \cite{Ma,Pet}.
For $n=2k$ even, and for any sequence $2\leq i_1<i_2<\cdots<i_k=n$, let $S_k(i_1,i_2,\ldots,i_k)$ denote the set of permutations in $\E_{2k}$ with peak set equals to $\{i_1,i_2,\ldots,i_k\}$, and let $s_k(i_1,i_2,\ldots,i_k)=|S_k(i_1,i_2,\ldots,i_k)|$. For $n=2k+1$ odd, and for any sequence $2\leq i_1<i_2<\cdots<i_k<i_{k+1}=n$, let $T_k(i_1,i_2,\ldots,i_{k+1})$ denote the set of permutations in $\E_{2k+1}$ with peak set equals to $\{i_1,i_2,\ldots,i_{k+1}\}$, and let $t_k(i_1,i_2,\ldots,i_{k+1})=|T_k(i_1,i_2,\ldots,i_{k+1})|$. Using induction on $n$, Strehl \cite{Str} derived the following enumeration formula for the numbers $s_k$ and $t_k$.
\begin{thm} [\cite{Str}]\label{thm}
\begin{align}\label{s}
s_k(i_1,i_2,\ldots,i_k)&= \prod_{1\leq j \leq k-1}(i_j-2j+1)^2, \\ \label{t}
t_k(i_1,i_2,\ldots,i_{k+1})&=\prod_{1\leq j \leq k}(i_j-2j+2)(i_j-2j+1).
\end{align}
\end{thm}
In this note we will give a combinatorial proof for the identities \eqref{s} and
\eqref{t}.

\section{Combinatorial Proof of Theorem \ref{thm}}

In order to give a combinatorial proof of Theorem \ref{thm}, we begin by deducing a formula for the number of matchings on $[2k]$ with a given closer set. Recall that a matching $\pi$ of $[2k]$ is a partition of the set $[2k]$ with the property that each block has exactly two elements. It can be represented as a graph with vertices $1, 2, \ldots, 2k$ drawn on a horizontal line in increasing order, where two vertices $i$ and $j$ are connected by an edge if and only
if $\{i,j\}$ (with $i < j$) is a block, and we say that $i$ is the
opener and $j$ is the closer of this edge. Since the graph is undirected, the edge can be denoted by $(i,j)$ or $(j,i)$ with no difference. The set of all the closers (resp. openers) of a matching is called its closer set (resp. opener set).

It is known that the number of matchings on $[2k]$ equals $(2k-1)!!$. If in addition, the $k$ closers among these $2k$ vertices are given, we have
\begin{lem} \label{pm}
The number of matchings on $[2k]$ with closer set $\{i_1<i_2<\cdots<i_k\}$ equals
\[
\prod_{1\leq j \leq k-1}(i_j-2j+1).
\]
\end{lem}
\pf Once the closer set is given, the opener set is accordingly known, and a matching on $[2k]$ can be constructed step by step by determining the $k$ edges with these closers $\{i_1,i_2,\cdots,i_k\}$ from left to right. For the closer $i_1$, there are $i_1-1$ openers before it with labels $1,2,\ldots,i_1-1$, so there are $i_1-1$ choices for the first closer to form an edge. Generally, for $2\leq j \leq k-1$, there are $j-1$ closers before the closer $i_j$, thus $j-1$ openers have been chosen so far by these closers, so the $j$-th closer has $i_j-1-2(j-1)$ openers to be chosen to generate a new edge. For the last closer $i_k$, there is only one opener left to form an edge with it. By considering all the $k$ closers, we see that there are $\prod_{1\leq j \leq k-1}(i_j-2j+1)$ possibilities to form a matching with this given closer set. \qed

As in \cite{Cor}, we can also represent a permutation of $[n]$ as a graph on $n$ vertices labeled $1, 2, \ldots , n$, with an edge from $i$ to $j$ if and only if $\sg_i = j$. Explicitly, put the $n$ vertices on a horizontal line, ordered from left to right by increasing label and we draw an edge from $i$ to $\sg_i$ above the line if $i\leq \sg_i$ and under
the line otherwise. Using this drawing, cycle up-down permutations having only even cycles correspond precisely to a pair of independent matchings whose vertices agree on openers and closers. By convention, we refer the matching with edges above the line as above matching and the matching with edges below the line as below matching. By this representation of permutations, we are now in the position to give a  combinatorial proof of Theorem \ref{thm}.

\noindent{\it Combinatorial Proof of \eqref{s}.}
Given an alternating permutation $\sg \in \E_{2k}$ with peak set $\{i_1,i_2,\ldots,i_k\}$, the reverse permutation $\sg^{r}$ is an up-down permutation on $[2k]$, and $\sg'=\tau(\sg^{r})$ is a cycle up-down permutation with only even cycles. Let $\G(\sg')$ be the corresponding graph of the permutation $\sg'$, it is easy to check that the closer set of the pair of matchings is $\{i_1,i_2,\ldots,i_k\}$. On the other hand, given a pair of matchings with closer set $\{i_1,i_2,\ldots,i_{k}\}$, we can recover an alternating permutation on $[2k]$ with peak set $\{i_1,i_2,\ldots,i_{k}\}$ by reversing the above procedure. By Lemma \ref{pm}, the number of above matchings and the number of below matchings with closer set $\{i_1,i_2,\ldots,i_{k}\}$ are both $\prod_{1\leq j \leq k-1}(i_j-2j+1)$, thus identity \eqref{s} follows. \qed

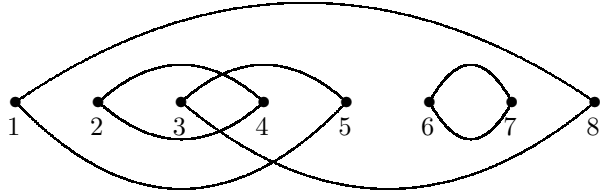
\begin{figure}[h,t]
\setlength{\unitlength}{1.1mm}
\begin{center}
\begin{picture}(80,20)

\multiput(2,6)(10,0){8}{\circle*{1.2}}

\put(1,2){\small 1}     \put(11,2){\small 2}    \put(21,2){\small 3}
\put(31,2){\small 4}    \put(41,2){\small 5}    \put(51,2){\small 6}
\put(61,2){\small 7}    \put(71,2){\small 8}

\qbezier[500](12,6)(22,15)(32,6) \qbezier[500](22,6)(32,15)(42,6)
\qbezier[500](52,6)(57,15)(62,6) \qbezier[500](2,6)(37,30)(72,6)

\qbezier[500](12,6)(22,-3)(32,6) \qbezier[500](52,6)(57,-3)(62,6)
\qbezier[500](2,6)(22,-15)(42,6) \qbezier[500](22,6)(47,-15)(72,6)
\end{picture}
\end{center}
\caption{The graph of the cycle up-down permutation (6,7)\,(2,4)\,(1,8,3,5).}\label{fig1}
\end{figure}

For example, if $\sg=5\,3\,8\,1\,4\,2\,7\,6 \in \E_8$, then $\sg^{r}=6\,7\,2\,4\,1\,8\,3\,5$ and $\tau(\sg^{r})=(6,7)\,(2,4)\,(1,8,3,5)$, the graph $\G(\tau(\sg^{r}))$ is depicted in Figure \ref{fig1}.

\noindent{\it Combinatorial Proof of \eqref{t}.}
Given an alternating permutation $\sg \in \E_{2k+1}$ with peak set $\{i_1,i_2,\ldots,i_{k+1}\}$, let $\overline{\sg}:=\sg\,0$ be a permutation on the set $\{0,1,2,\ldots 2k,2k+1\}$ obtained by appending a ``0" after the last element of $\sg$. Since $\overline{\sg}_{2k}=\sg_{2k}<\overline{\sg}_{2k+1}=\sg_{2k+1}>0=\overline{\sg}_{2k+2}$, and $\overline{\sg}_{i}=\sg_i$ for $i\leq 2k+1$, we can view $\overline{\sg}$ as an alternating permutation of $\E_{2k+2}$ and the reverse permutation $\overline{\sg}^{r}$ is an up-down permutation on $[2k+2]$ with the first element being $0$. From the position of $0$, we see that $0$ is the unique left to right mimimum of $\overline{\sg}$, thus $\overline{\sg}'=\tau(\overline{\sg}^{r})$ is a cycle up-down permutation with only one even cycle.

Let $\G(\overline{\sg}')$ be the corresponding graph of the permutation $\overline{\sg}'$, then the closer set of the pair of matchings is $\{i_1,i_2,\ldots,i_{k+1}\}$. Since there are only one cycle in $\overline{\sg}'$, the above matching and the below matching are not independent now. It requires that in the graph $\G(\overline{\sg}')$, the edges not containing the last closer $i_{k+1}$ can not form a closed circle. That is to say, for any closer $i_j, j\leq k$, there exists no opener $a$ such that $(a,b_1),(b_1,b_2),\cdots, (b_m,i_j)(i_j,a)$ are the edges of $\G(\overline{\sg}')$. On the other hand, given a pair of matchings with closer set $\{i_1,i_2,\ldots,i_{k+1}\}$ and satisfy the above condition, we can also reverse the above procedure to get an alternating permutation on $[2k+1]$ with peak set $\{i_1,i_2,\ldots,i_{k+1}\}$.

By the same analysis as Lemma \ref{pm}, the number of above matchings on $\{0,1,2,\ldots 2k,2k+1\}$ with closer set $\{i_1,i_2,\ldots,i_{k+1}\}$ equals $\prod_{1\leq j \leq k}(i_j-2j+2)$ since there is an extra opener $0$. After the above matching is determined, we next construct the edges of the below matching. For $1\leq j \leq k$, there are $i_j-1-2(j-1)+1$ openers before the $j$-th closer $i_j$. Among these openers, the opener $a$ such that $(i_j,a),(a,b_1),(b_1,b_2),\cdots , (b_m,i_j)$ are the edges of the graph have been constructed so far can not be chosen, otherwise there will exist at least two cycles in $\overline{\sg}'$. Hence the $j$-th closer of the below matching has $i_j-2j+1$ choices for its opener to generate an edge. Combining the number of possible ways of constructing this pair of matchings with the same given closer set leads to the identity \eqref{t} immediately.   \qed

\begin{figure}[h,t]
\setlength{\unitlength}{1.1mm}
\begin{center}
\begin{picture}(100,20)
\multiput(2,6)(10,0){10}{\circle*{1.2}}
\put(1,2){\small 0}     \put(11,2){\small 1}    \put(21,2){\small 2}
\put(31,2){\small 3}    \put(41,2){\small 4}    \put(51,2){\small 5}
\put(61,2){\small 6}    \put(71,2){\small 7}    \put(81,2){\small8}
\put(91,2){\small 9}
\qbezier[1000](2,6)(27,30)(52,6) \qbezier[1000](22,6)(59,30)(92,6)
\qbezier[1000](12,6)(27,20)(42,6) \qbezier[1000](32,6)(52,20)(72,6)
\qbezier[1000](62,6)(72,20)(82,6)
\qbezier[1000](32,6)(37,-2)(42,6)  \qbezier[1000](22,6)(37,-10)(52,6)
\qbezier[1000](62,6)(67,-2)(72,6) \qbezier[1000](2,6)(42,-20)(82,6) \qbezier[1000](12,6)(52,-20)(92,6)
\end{picture}
\end{center}
\caption{The graph of the cycle up-down permutation (0,5,2,9,1,4,3,7,6,8).}\label{fig2}
\end{figure}
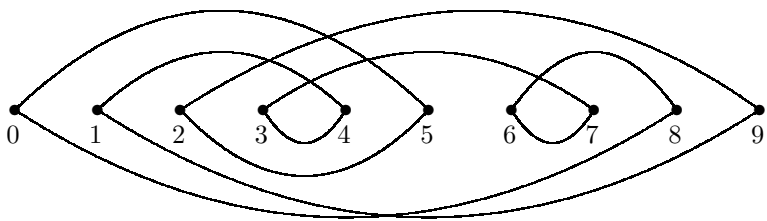

Let us illustrate the argument with an example. For $\sg=8\,6\,7\,3\,4\,1\,9\,2\,5 \in \E_9$, then $\overline{\sg}=8\,6\,7\,3\,4\,1\,9\,2\,5\,0\in \E_{10}$,  $\overline{\sg}^{r}=0\,5\,2\,9\,1\,4\,3\,7\,6\,8$ and $\tau(\overline{\sg}^{r})=(0,5,2,9,1,4,3,7,6,8)$, see Figure \ref{fig2} for its corresponding graph.

\end{document}